\documentclass[11pt,leqno]{amsart}
\usepackage{amssymb}
\newtheorem{lemma}{Lemma}
\newtheorem{corollary}{Corollary}
\newtheorem*{combinatoriallemma}{Combinatorial Lemma}
\theoremstyle{definition}
\newtheorem{definition}{Definition}
\newtheorem{observation}{Observation}
\def\st{such that}
\def\wrt{with respect to}
\def\ie{{\em i.e.,}} 
\def\NN{\mathbb{N}}
\def\sC{\mathcal{C}}
\def\sP{\mathcal{P}}
\def\sS{\mathcal{S}}
\def\sV{\mathcal{V}}

\begin{document}
\title[Aumann's Theorem]{On Aumann's Theorem that the sphere does not 
admit a mean}
\author[F. J. Trigos-Arrieta]{F. Javier Trigos-Arrieta}
\address{Department of Mathematics, California State University, 
Bakersfield, 9011 Stockdale Highway, Bakersfield, 
California 93311-1099, USA}
\email{jtrigos@csubak.edu}
\author[M. Turza\'nski]{Marian Turza\'nski}
\address{Department of Mathematics, Cardinal Stefan Wyszy\'nski University,
ul. Dewajtis 5, 01-815 Warsaw, Poland} 
\email{mtturz@ux2.math.us.edu.pl}
\date{June 2003.  Work in Progress}
\begin{abstract} 
We prove that the circle $S_1$ does not have a 2-mean, i.e.,
$S_1\times S_1$ cannot have a retraction $r$ onto its diagonal
with $r(x,y)=r(y,x)$, whenever $x,y \in S_1$. Our proof is
combinatorial and topological rather than analytical.
\end{abstract} 
\maketitle

\section{Introduction}

{\sc Aumann} and {\sc Caratheodory} \cite{AuCa34}, \cite{Au35} and
\cite{Au43} were among the pioneers who first considered the question
about the structure of spaces for which the topological product $X^n$ has
a symmetric retraction onto its diagonal, \ie \ {\em a $n$-mean}. They
studied such objects in the complex plane and in the Euclidean $n$-space
using analytical tools. For example {\sc Aumann} in \cite{Au43} proved
that the $n$-dimensional sphere does not have a mean. For more information
about means see \cite{Hilton}.  The aim of this note is to prove that the
circle $S_1$ does not have a 2-mean, using only combinatorial and
topological tools.  For this purpose we use a method comparable to the one
used in \cite{kuS}.  It is interesting to notice that this method has been
used (in dimension 2)  to prove, among some other results, the {\sc
Brouwer} fixed point theorem and the special hexagonal chessboard theorem
(see {\sc Gale} \cite{Gal}, who, as far as we know, introduced the
method), and the {\sc Borsuk-Ulam} antipodal theorem (see \cite{KuTu}).  
We point out that in the case of the {\sc Brouwer} fixed point theorem,
the combinatorial proof in \cite{KKM} is based on {\sc Sperner's} Lemma
\cite{Sp} but in the case of the {\sc Borsuk-Ulam} antipodal theorem
\cite{Bor}, for the combinatorial proof {\sc Tucker's} Lemma
\cite{Tucker45} is used (in this case {\sc Sperner's} Lemma is not enough)
({\sc Ky Fan} \cite{Fan52} extended {\sc Tucker}'s result to arbitrary
$n$). For more information about fixed point theory see \cite{D-G}. Here
we have (in dimension 2) one universal combinatorial lemma (see next
section). We wonder if it is possible to generalize this method to
arbitrary $n$.

\section{Combinatorial part}

Let  us fix a natural number $k > 1$ and let
$$
Z_k = \left\{ \frac{i}{k}: i \in \{0,..., k \} \right\}
$$
and denote by 
$$
D^{2}(k) 
= (Z_k \times  Z_k) 
= \left\{0,\frac{1}{k},..., \frac{k-1}{k},1\right\}^2;
$$ 
$D^{2}(k)$ is called {\em a combinatorial square.}

\begin{definition}
Denote by ${\mathbf e}_0 = ( \frac{1}{k}, 0), {\mathbf e}_1 =  ( 0, 
\frac{1}{k})$  the basic vectors of length $\frac{1}{k}$. An ordered set 
$z = [z_0, z_1, z_2]$ is said to be a {\em simplex} if and only if
$$z_1 = z_0 + \mathbf{e}_i, z_2 = z_1 + \mathbf{e}_{1-i}
\quad\text{where $i \in \{ 0, 1 \}$.}$$ 
Any subset $[z_0, z_1 ], [z_1, z_2 ]$ and $[z_2, z_0] \subset z$ is 
said to be a face of the simplex $z$.
\end{definition}  

\begin{figure}[h]
\setlength{\unitlength}{1cm}
\begin{picture}(3,2)
\thicklines
\put(.8,0){$\rightarrow$} %1
\put(.4,.1){\line(1,0){.7}} %1
\put(1.32,.3){\line(0,1){.7}} %1
\put(.3,.25){\line(1,1){.75}} %1
\put(1.2,0){$z_1$} %1
\put(1.23,.8){$\uparrow$} %1
\put(0,0){$z_0$} %1
\put(1.1,1.2){$z_2$} %1
\put(.23,.25){$\swarrow$} %1
\put(0.6,-.3){$\mathbf{e}_0$} %1
\put(1.6,.5){$\mathbf{e}_1$} %1
\put(2.35,1.3){\line(1,0){.65}} %2
\put(2.15,.3){\line(0,1){.7}} %2
\put(2.3,.25){\line(1,1){.75}} %2
\put(2,0){$z_0$} %2
\put(2,1.2){$z_1$} %2
\put(3.1,1.2){$z_2$} %2
\put(2.7,1.2){$\rightarrow$} %2
\put(2.05,.8){$\uparrow$} %2
\put(2.23,.25){$\swarrow$} %1
\put(2.6,1.5){$\mathbf{e}_0$} %1
\end{picture}
\caption{}
\end{figure}

\begin{observation}
Any face of a simplex $z$ contained in  $D^{2}(k)$ 
is a  face of exactly one or two simplexes from $D^{2}(k)$, depending on 
whether or not it lies on the boundary of $D^{2}(k)$.
\end{observation}

\begin{definition}
Let ${\mathcal{P}}(k)$ be the family of all simplexes in $D^{2}(k)$ and 
let ${\mathcal{V}}(k)$ be the set of all vertices of the simplexes from
${\mathcal{P}}(k)$. {\em A coloring} of ${\mathcal{P}}(k)$ is any 
function $f:{ \mathcal{V}}(k) \longrightarrow \{1, -1\}$, and any face $s$ 
of any simplex $z$ is called {\em an $f$-gate} (or simply a gate if there 
is no ambiguity of what $f$ is) if $f[s] = \{1, -1\}$.
\end{definition}

\begin{observation}
Let $w$ be a simplex, $\mathcal{W}$ be the set of vertices of $w$
and $f: {\mathcal{W}} \longrightarrow  \{1, -1\}$ be a function. Then 
$w$ has an even number of gates.
\end{observation}

\begin{definition}
If $f:{ \mathcal{V}}(k) \longrightarrow \{0, 1\}$ is a function, two 
simplexes $w$ and $v$ from ${\mathcal {P}}(k) $ are in the relation $\sim$ 
if $w \cap v$ is a gate. A subset ${\mathcal{S}} \subset {\mathcal{P}}(k)$ 
is called a chain in ${\mathcal{P}}(k)$ if ${\mathcal{S}} = \{w_0, 
w_1,..., w_n\}$ and for each $i\in\{0,..., n-1\},  w_{i} \sim w_{i+1}$.
\end{definition}

\begin{observation}
For each chain $\{v_1,...,v_n \}\subset {\mathcal{P}}(k)$  there exists no 
more than one $v \in \sP(k)$ and one  $w \in \sP(k)$ such that  
$\{v_1,...,v_n,v \}$ and $\{w,v_1,...,v_n \}$ are chains. 
Also, if ${\mathcal{S}}_1$ and ${\mathcal{S}}_2$ are maximal chains in 
${\mathcal{P}}(k)$, then either 
${\mathcal{S}}_1 \cap {\mathcal{S}}_2 = \emptyset$ or  
${\mathcal{S}}_1 = {\mathcal{S}}_2$.
\end{observation}

Let $a$ and $b$ be two different elements of $D^2(k)$.  
Consider the rectangle $R$ with $a$ and $b$ as opposites vertices and
right-hand-orient its boundary. By $\overline{ab}$ we mean the part of 
the boundary that goes from $a$ to $b$. We define similarly $\overline{ba}$.
The boundary of $R$ is denoted by $\partial R$.

\begin{lemma}
No maximal chain $\sS \subseteq \sP(k)$ ever finishes at a gate of an 
interior simplex, \ie \ a simplex disjoint from $\partial R$.
\end{lemma}

\begin{proof}
It consists to show that if a simplex 
$S_1$ in $\sS$ is disjoint with $\partial R$, there is always another 
simplex $S_2$ with a common gate (Observation 2). Thus the only 
possibility for $\sS$ to stop is at $\partial R$. 
There are twelve possible (simplex, flow of the chain (if directed)) 
combinations of the simplex to be considered,
each of them with two possible outcomes. We picture some of them with the
following in mind: arrows mean flow, thick lines are NOT gates and
thin lines are gates:

\begin{figure}[h]
\setlength{\unitlength}{1cm}
\begin{picture}(8,2)
\thicklines
\put(.35,1.1){\line(1,0){.65}} %1
\put(2.35,1.1){\line(1,0){.65}} %2
\put(3.15,0.29){\line(0,1){.65}} %2
\put(4.35,1.1){\line(1,0){.65}} %3
\put(4.35,.1){\line(1,0){.6}} %3
\thinlines
\put(.15,.3){\line(0,1){.65}} %1
\put(2.15,.3){\line(0,1){.65}} %2
\put(2.35,.1){\line(1,0){.6}} %2
\put(4.15,.3){\line(0,1){.65}} %3
\put(.3,.25){\line(1,1){.75}} %1
\put(2.3,.25){\line(1,1){.75}} %2
\put(4.3,.25){\line(1,1){.75}} %3
\put(5.15,0.29){\line(0,1){.65}} %3
\put(0,0){$\oplus$} %1
\put(2,0){$\oplus$} %2
\put(4,0){$\oplus$} %3
\put(0,1){$\ominus$} %1
\put(2,1){$\ominus$} %2
\put(4,1){$\ominus$} %3
\put(1,1){$\ominus$} %1
\put(3,1){$\ominus$} %2
\put(5,0){$\oplus$} %3
\put(5,1){$\ominus$} %3
\put(3,0){$\ominus$} %2
\put(-.05,.5){$\rightarrow$} %1
\put(1.95,.5){$\rightarrow$} %2
\put(3.95,.5){$\rightarrow$} %3
\put(4.95,.5){$\rightarrow$} %3
\put(.55,.5){$\searrow$} %1
\put(2.55,.5){$\searrow$} %2
\put(4.55,.5){$\searrow$} %3
\put(2.55,-.1){$\downarrow$} %2
\put(.55,1.5){A} %1
\put(2.5,1.5){A1} %1
\put(4.5,1.5){A2} %1
\end{picture}

\setlength{\unitlength}{1cm}
\begin{picture}(8,2)
\thicklines
\put(.35,1.1){\line(1,0){.65}} %1
\put(3.35,1.1){\line(1,0){.65}} %2
\put(6.35,1.1){\line(1,0){.65}} %3
\put(2.2,0.3){\line(1,1){.75}} %2
\put(5.35,.15){\line(1,0){.6}} %3
\thinlines
\put(.15,.3){\line(0,1){.65}} %1
\put(3.15,.3){\line(0,1){.65}} %2
\put(6.15,.3){\line(0,1){.65}} %3
\put(.3,.25){\line(1,1){.75}} %1
\put(3.3,.25){\line(1,1){.75}} %2
\put(6.3,.25){\line(1,1){.75}} %3
\put(5.2,0.3){\line(1,1){.75}} %3
\put(2.35,.15){\line(1,0){.6}} %2
\put(0,0){$\oplus$} %1
\put(3,0){$\oplus$} %2
\put(6,0){$\oplus$} %3
\put(2,0){$\ominus$} %2
\put(5,0){$\oplus$} %3
\put(0,1){$\ominus$} %1
\put(3,1){$\ominus$} %2
\put(6,1){$\ominus$} %3
\put(1,1){$\ominus$} %1
\put(4,1){$\ominus$} %2
\put(7,1){$\ominus$} %3
\put(-.05,.5){$\leftarrow$} %1
\put(2.95,.5){$\leftarrow$} %2
\put(5.95,.5){$\leftarrow$} %3
\put(5.25,.5){$\leftarrow$} %3
\put(.55,.5){$\nwarrow$} %1
\put(3.55,.5){$\nwarrow$} %2
\put(6.55,.5){$\nwarrow$} %3
\put(2.55,0){$\downarrow$} %2
\put(.55,1.5){B} %1
\put(3.55,1.5){B1} %2
\put(6.55,1.5){B2} %3
\end{picture}
\caption{}
\end{figure}

\end{proof}

\begin{corollary}
Any maximal chain $\sS \subseteq \sP(k)$ beginning at $\partial R$ must 
finish at $\partial R$.
\end{corollary}

\begin{combinatoriallemma}
Let $\sP(k)$ be the set of simplexes of
$D^2(k)$ and $f:\sV(k)\to \{-1,1\}$ be a coloring of $\sV(k)$. If 
$a$ and $b$ belong to $\sV(k)$, and
$f(b)=-f(a)$, then there exists a chain $\sS \subseteq \sP(k)$ such
that $\sS \cap \overline{ab} \neq \emptyset \neq \sS \cap \overline{ba}$.
\end{combinatoriallemma}

This result was proved originally in \cite{Tur1}. 
Here we present a different argument.

\begin{proof}
We first define two equivalence relations on $\sV(k)$:

If $u,v \in D^2(k) \cap R$, we will say that $u\approx v$ if $u=v$ or 
if there
are vertices $u=x_0, x_1,...,x_{n-1},x_n=v$ in $D^2(k) \cap R$ \st \
$[x_i,x_{i+1}]$ is a face of a simplex ($i=0,...,n-1$) and 
$f(x_i)=f(x_{i+1})$. Clearly $\approx$ is an equivalence relation on
$D^2(k) \cap R$.

Let $\sS \subseteq \sP(k)$ be a maximal chain beginning 
at the boundary of $R$. If $u,v \in D^2(k) \cap R$, we will say that 
$u\simeq v$ if $u=v$ or 
if there are vertices $u=x_0, x_1,...,x_{n-1},x_n=v$ in $D^2(k) \cap R$ 
with $[x_i,x_{i+1}]$ being a face of a simplex ($i\in \{0,...,n-1\}$) and
no $[x_i,x_{i+1}]$ is a gate belonging to a simplex belonging to $\sS$. 
Clearly $\simeq$ is an equivalence relation on $D^2(k) \cap R$ as well.

Let $\sC$ be the $\approx$-component of $b$. Walking from $a$ to $b$ let 
$x$ be the vertex on $\overline{ab}$ found right  before $\sC \cap 
\overline{ab}$, and $y$ be the vertex on $\overline{ab}$ right after $x$. 
Then $y \in \sC$ and $f(x)=f(a)$. Thus $[x,y]$ is a gate.
Let $\sS$ be the unique maximal chain to which the simplex containing
$[x,y]$ belongs to (Observation 3). By Corollary 1, $\sS$ ends on 
$\partial R$. By the choice of $x$ and $y$ and since points in $\sC$ are 
all $\simeq$-equivalent, $\sS$ must end on  $\overline{ba}$, as required.
\end{proof}
 
\section{Topological Part}

We borrow the following from \cite{kuS}.

\begin{definition}
If $\{A_m: m \in \NN \}$ is a sequence of subsets of a compact metric 
space $X$, we define its {\em upper limit} $Ls\{A_n: n \in \NN\}$ as the
set of points $x \in X$ \st \ there is an infinite $M \subseteq \NN$  
such that for every $m \in M$ there is $x_m \in A_M$ with $x_m \to x$.
\end{definition}

In the paper \cite{kuS} the following result has been proved. See also
\cite{Kur2} (5.47.6).

\begin{lemma}
Let $\{A_m: m \in \NN \}$ be a sequence of connected subsets of a compact 
metric space $X$ such that some sequence $\{a_n: n \in \NN\}$ of points 
$a_n \in A_n$ is converging in $X$. Then the set $Ls\{A_n: n \in N\}$ is 
compact and connected.
\end{lemma}

\section{Main Result}

In this section we prove the result mentioned in the abstract.

Let $X$ be a space, and denote by $\Delta(X^2):=\{(x,x):x\in X\}$. 
Obviously $\Delta(X^2)$ is homeomorphic to $X$. Identify $S_1$ with 
$I:=[0,1]$ and $0=1$.

Suppose that there exists a symmetric retraction $r$ from $S_1\times S_1$ 
onto its diagonal $\Delta(S_1^2)$, \ie \ a continuous map
$r:S_1\times S_1 \to \Delta(S_1^2)$ satisfying:\\
a) $r(x,y)=r(y,x)$ for each $x$ and $y$ from $S_1$, and \\
b) $r(x,x)=(x,x)$. \\
We call $r$ a {\em 2-mean,} and say that $S_1$ has a 2-mean.
 
To prove that the existence of $r:S_1^2\to \Delta(S_1^2)$ with properties
(a-b) is impossible, we consider two cases:

(1) Assume that $r[(I\times \{1\}) \cup (\{0\}\times I)]\neq \{(0,0)\}$.

Notice that if we consider $I^2$ instead of $S_1^2$, and $r:I^2\to 
\Delta(I^2)$ rather than $r:S_1^2\to \Delta(S_1^2)$, then $r$ has the 
following additional properties: \\
c) $r(0,0)=(0,0),\:r(1,1)=(1,1)$,\\
d) $r(0,x)=r(1,x),\:r(x,0)=r(x,1).$

For illustrative purposes, we call $\{0\}\times I$ 
$:=$ ``left'', $\{1\}\times I$ $:=$ ``right'', $I\times \{0\}:=$ 
``bottom'' and
$ I \times \{1\}:=$ ``top''. The assumption we are assuming reads now
$r[(I\times \{1\}) \cup (\{0\}\times I)]\neq \{(0,0),(1,1)\}$. (c) and the
Intermediate Value Theorem imply that 
$r[(I\times \{1\}) \cup (\{0\}\times I)]=\Delta(I^2)$.
Fix $k \in \NN$, and if 
$p:I\times I \to I$ denotes the projection on the first coordinate, 
define the coloring $f:V(k)\to \{\pm1\}$ as follows:
$$f(i/k,j/k)=\left\{\begin{array}{rc}
		-1&  \mbox{if}\:\: \cos(2\pi p(r(i/k,j/k)))\leq 0,\\ 
		1&   \mbox{otherwise.}
		\end{array} \right.$$
This coloring is symmetric \wrt \
$\Delta(I^2)$ and each side of the square has exactly the same number of 
gates: The gates at the left and right sides are at the same vertical 
positions, and those at the bottom and top sides are at the same 
horizontal positions, respectively. 

Considering once again $r:S_1^2\to \Delta(S_1^2)$, we identify the points 
$(0,i/k)$ with $(1,i/k)$ and $(i/k,0)$ with $(i/k,1)$ ($i=0,...,k$). 
Walking to the right of $(0,0)$, one finds the first gate $g_b^1$ 
($b$, $l$, $r$ and $t$ 
stand for ``bottom'',''left'', ``right'' and ``top'') 
on $I \times \{0\}$ which gives place to a chain $\sS_k$ ``going'' on top
of the $\approx$-component $A$ of $(0,0)$ (the relation $\approx$ was 
defined in the proof of the Combinatorial Lemma). By the case we are 
dealing with, $\sS_k$ intersects $\{0\}\times I$
in the last gate $g_l^\infty$ from top to bottom.
By the identification of $(0,i/k)$ with $(1,i/k)$ $\sS_k$ reappears through 
the first gate $g_r^1$ in $\{1\}\times I$ from bottom to top, and thus 
$\sS_k$ ``goes'' above the $\approx$-component $B$ of $(1,0)$. 
 
$\sS_k$ intersects $I \times \{0\}$ in the last gate $g_b^\infty$ 
going from left to right. By the identification of $(0,i/k)$ with 
$(1,i/k)$ $\sS_k$ reappears through 
the first gate $g_t^1$ of the top side from right to left, and thus 
$\sS_k$ ``goes'' under
the $\approx$-component $C$ of $(1,1)$. Again $\sS_k$ intersects 
$I \times \{1\}$ in the last gate $g_r^\infty$ of the right side
going from bottom to top,
thus $\sS_k$ reappears on the first gate $g_l^1$ of the left side
going from top to bottom,
going under the $\approx$-component $D$ of $(0,1)$ and intersecting 
the last gate $g_t^\infty$ of the top side
from right to left, reappearing on $g_b^1$
and beginning the whole cycle once again. 
The union of the simplexes from the chain $\sS_k$ is a connected set
for each natural number $k$.

According to Lemma 2 the upper limit $C=Ls\{\sS_k: k\in \NN\}$ is 
connected, and we have that $C\subset r^{-1}(p^{-1}(\cos^{-1}(0)))$, thus 
$r$ maps the continuum $C$ onto two points in $\Delta(S_1^2)$; a 
contradiction.

This concludes the proof in case (1).

(2) Assume that $r[(I\times \{1\}) \cup (\{0\}\times I)]= \{(0,0)\}$.

This would mean that $r[\partial I^2]=\{(0,0)\}$, and thus would imply
that any copy $S$ of $S_1$
in the sphere $S_2$ is a retract: The sphere $S_2$ is the image of
$S_1 \times S_1$ by identifying the left (and right) and 
bottom (and top) sides of $S_1 \times S_1$. 
Since $S_2$ equals the union of two copies of the unit disk,
sharing the same boundary, this is impossible by the following corollary 
to the Combinatorial Lemma:

\begin{corollary}[Borsuk's non-retraction theorem]
$S_1$ is not a retract of the unit disk.
\end{corollary}

\begin{proof}
Identify the disk with the square $I^2$, and $S_1$ with  its boundary 
$\partial I^2$.
If $k \in \NN$ consider $D^2(k)$ and color it according to what points
get mapped to the bottom and left sides, and to the top and right sides.
There are only two gates in $\partial I^2$. By corollary 1 there is one and
only one chain connecting these two gates. Then we proceed similarly
as in the end of case (1).
\end{proof}

\bibliographystyle{amsplain}
\makeatletter 
\renewcommand{\@biblabel}[1]{\hfill#1.}
\makeatother

\end{document}